\documentclass{ifacconf}

\usepackage{graphicx}      
\usepackage{natbib}        
\usepackage{amsmath,amssymb,amsfonts}
\newtheorem{assumption}{$Assumption$}
\newtheorem{lemma}{$Lemma$}
\newtheorem{remark}{$Remark$}
\newtheorem{definition}{$Definition$}
\usepackage{algorithm}
\usepackage{algorithmic}

\begin{document}
\begin{frontmatter}

\title{Towards an $O(\frac{1}{t})$ convergence rate for distributed dual averaging}


\author[First]{Changxin Liu} 
\author[Second]{Huiping Li} 
\author[First]{Yang Shi}

\address[First]{Department of Mechanical Engineering, University of Victoria, Victoria, B.C., Canada (e-mail: chxliu@uvic.ca, yshi@uvic.ca)}
\address[Second]{School of Marine Science and Technology, Northwestern Polytechnical University, Xi'an, P. R. China (e-mail: lihuiping@nwpu.edu.cn)}

\begin{abstract}                
Recently, distributed dual averaging has received increasing attention due to its superiority in handling constraints and dynamic networks in multiagent optimization. However, all distributed dual averaging methods reported so far considered nonsmooth problems and have a convergence rate of $O(\frac{1}{\sqrt{t}})$. To achieve an improved convergence guarantee for smooth problems, this work proposes a second-order consensus scheme that assists each agent to locally track the global dual variable more accurately. 
This new scheme in conjunction with smoothness of the objective ensures that the accumulation of consensus error over time caused by incomplete global information is bounded from above.
Then, a rigorous investigation of dual averaging with inexact gradient oracles is carried out to compensate the consensus error and achieve an $O(\frac{1}{t})$ convergence rate. The proposed method is examined in a large-scale LASSO problem.
\end{abstract}


\end{frontmatter}

\section{Introduction}
%
%
We consider the problem where a team of agents connected via a network manage to optimize the sum of their local interests while respecting certain common constraints. This problem is referred to as distributed optimization, and has been extensively investigated in recent years mainly due to its broad applications. For example, distributed machine learning, formation control of autonomous vehicles, and sensor fusion can be cast as optimization problems of this type. For a recent overview of distributed optimization, please refer to \cite{nedic2018network}.

In such a framework, each agent does not have full knowledge of the objective function, therefore it has to communicate with neighbors to estimate the global information, e.g.,  global gradient or/and mean value of local variables, during the course of optimizer seeking to achieve distributed optimization. 
Regarding the estimation process, the algorithms in \cite{nedic2010constrained,yuan2016convergence,liu2019distributed} directly seek consensus over local variables based on a doubly stochastic weight matrix, while the optimizer seeking process is guided by the gradient of the local objective. However, due to the fact that local gradients evaluated at global minimizer are not necessarily zero, the two forces caused by consensus and local gradient flows are conflicting with each other, preventing exact optimization when a constant stepsize is used, that is, there always exists a gap between the accumulation point and global minimum. It is worth mentioning that, by using a decaying stepsize, exact optimization may be obtained with, however, a slow $O(\frac{1}{\sqrt{t}})$ convergence rate where $t$ is the time counter. 
This issue can be solved by an additional estimation process for the global gradient by using the dynamic average consensus scheme in \cite{zhu2010discrete}. It is shown in \cite{varagnolo2015newton,qu2017harnessing} that for unconstrained smooth optimization the algorithm steered by the approximated global gradient obtains exact minimization with an $O(\frac{1}{{t}})$ rate.

In the methods mentioned above, local estimates about the minimizer are directly generated in the feasible set (in the case of constrained optimization) that is contained in the primal space of variables. There are also some schemes available in the literature where the minimizer seeking process imitates dual methods, e.g., mirror descent in \cite{shahrampour2017distributed} and dual averaging in \cite{duchi2011dual}. The concept of dual methods was coined by \cite{nemirovsky1983problem}, where a dual model of the objective is updated and a prox-function establishes a mapping from the dual space to the primal to shrink the error bound in primal methods.
For example, \cite{duchi2011dual} designed a distributed dual averaging algorithm where the global dual variable is gradually learned by a consensus scheme, and demonstrated that minimizing the approximate dual model of the global objective helps bypass the difficulty caused by projection in distributed primal methods. Recent work in \cite{liu2018distributed} introduced another averaging step to standard distributed dual averaging to reap a non-ergodic convergence property, which helps deal with distributed optimization problems with coupled constraints. For problems defined over time-varying and unbalanced networks, a distributed dual averaging method with the push-sum technique was reported in \cite{liang2019dual}.
 
Although distributed dual methods in the literature have demonstrated advantages
over their primal counterparts in terms of constraint handling, convergence
rate, and analysis complexity, all the results reported so far focused only on nonsmooth
optimization and have a convergence rate of $O(\frac{1}{\sqrt{t}})$. \emph {Considering this, a question naturally arises: If the objective functions
exhibit some desired properties, e.g., smoothness, is it possible to accelerate the convergence
rate of distributed dual averaging to $O(\frac{1}{t})$?}
This work provides affirmative answer to this question. This is made admissible by a new second-order consensus scheme that assists each agent to locally track the global dual variable more accurately. With the new dual estimate, the accumulation of error over time between local primal variables and their mean is proved to admit an upper bound. This together with a rigorous investigation of averaged primal variables yields an accelerated convergence rate.

\emph{Notation}: $\mathbb{R}$ represents the set of real numbers and $\mathbb{R}^m$ the $m$-dimensional Euclidean space. 
In this space, we let $\lVert \cdot \rVert_p$ denote the $l_p$-norm operator, and without specifying $p$, it stands for the Euclidean norm. 
We denote by
$0_m\in\mathbb{R}^m$
and $\mathbf{1}$
 the $m$-dimensional vector of all zeros and all ones, respectively.
Given a matrix $P \in\mathbb{R}^{m\times m}$, its spectral radius and singular values are denoted by $\rho(P)$ and $\sigma_1(P)\geq \sigma_2(P)\geq \cdots \geq \sigma_m(P) \geq 0$, respectively.



%
%
%
%
%
%
\section{Problem Statement and Preliminaries}

\subsection {Problem Statement}

Formally, the optimization problem is given by
\begin{equation} \label{constrained_optimization}
\min_{x\in\mathcal{X}} f(x) =  \sum_{i=1}^{n}f_i(x)
\end{equation}
where $x\in\mathbb{R}^m$ denotes the global decision variable, $f_i: \mathbb{R}^m\rightarrow\mathbb{R}, i\in\mathbb{N}_{[1,n]}$ represents the local objective function that is privately known by agent $i$, and $\mathcal{X}\subseteq\mathbb{R}^m$ stands for the common constraint set. Throughout the paper, we denote one of the minimizers by $x^*$. For \eqref{constrained_optimization}, we make the following standard assumption.

\begin{assumption}\label{lipschitzassumption}
	Each $f_i(x), i\in\mathbb{N}_{[1,n]}$ is convex and has Lipschitz continuous gradient with parameter $L$, i.e.,
	\begin{equation*}
	\|\triangledown f_i(x)- \triangledown f_i(y)\|\leq L\|x-y \|, \forall x,y\in\mathcal{X}.
	\end{equation*}
	The common constraint set $\mathcal{X}$ is convex and closed, and contains the origin.
\end{assumption}

We use an undirected graph $\mathcal{G}=\{\mathcal{V},\mathcal{E}\}$ to describe the communication pattern between agents, where $\mathcal{V}=\{1,\cdots,n\}$ denotes the set of $n$ agents and $\mathcal{E}\subseteq \mathcal{V}\times\mathcal {V}$ represents the set of channels that connect agents, that is, the pair $(i,j)\in \mathcal{E}$ for $i,j\in \mathcal{V}$ indicates that there exists a link between node $i$ and $j$. The set of $i$'s neighbors is denoted by $ \mathcal{N}_i=\{j\in \mathcal{V}|(j,i)\in \mathcal{E} \}$. The graph is assumed to be fixed and connected in the following.

\begin{assumption}\label{graphconnected}
	The communication graph $\mathcal{G}=(\mathcal{V},\mathcal{E})$ is fixed and connected.
\end{assumption}

Based on Assumption \ref{graphconnected}, a proper weight matrix $P=[p_{ij}]$ can be constructed.
In particular, a positive weight $p_{ij}$ is assigned to each communication link $(i,j)\in\mathcal{E}$; for other $(i,j)$ pairs, zero weight is considered. Moreover, the weight matrix satisfies the following assumption.
\begin{assumption}\label{assumption_weight_matrix}
	1) $P$ has a strictly positive diagonal, i.e., $p_{ii}>0$;
	2) $P$ is doubly stochastic, i.e., $P\mathbf{1}=\mathbf{1}$ and $\mathbf{1}^{\mathrm{T}}P=\mathbf{1}^{\mathrm{T}}$.
\end{assumption}

Without loss of generality, we will assume $m=1$ for ease of notation in the remaining sections, i.e., $\mathbf{1}\otimes I_m=\mathbf{1}$, $P\otimes I_m=P$.

\subsection {Preliminaries}

\begin{definition}\label{prox_function}
	A function $d: \mathcal{X}\rightarrow\mathbb{R}$ is called a prox-function if 1) $d(x)\geq 0, \forall x\in\mathcal{X}$ and $d(0_m)=0$;
	2) $d(x)$ is differentiable and $1$-strongly convex on $\mathcal{X}$, i.e.,
	\begin{equation*}
	d(y)\geq d(x)+\langle\triangledown d(x), y-x \rangle+\frac{1}{2}\lVert y-x\rVert^2, \forall x,y\in\mathcal{X}.
	\end{equation*}
\end{definition}

\begin{definition}\label{bregman_divergence}
	For $x,y\in\mathcal{X}$, the Bregman divergence induced by a prox-function $d$ is defined as
	\begin{equation*}
	D_{d}(x,y)=d(x)-d(y)-\langle \triangledown d(y), x-y \rangle.
	\end{equation*}
\end{definition} 

%
%

\section{Algorithm Development}
\subsection {Centralized Dual Averaging}
This subsection introduces the centralized dual averaging (${\rm CDA}$) \cite{nesterov2009primal}. 
${\rm CDA}$ generates sequences of the estimates  
about the minimizer ($\{x_t\}_{t\geq0}$) and the dual variable ($\{\sum_{k=0}^{t}\triangledown f(x_k)\}_{t\geq0}$) according to the following rule:
\begin{equation}\label{cda}
	{x}_{t+1}= \arg \min_{x\in\mathcal{X}}\big\{a_{t}\sum_{k=0}^{t}\langle \triangledown f(x_{k}),x \rangle+d(x)\big\}
\end{equation}
where $\{a_t\}_{t\geq 0}$ is a sequence of positive control parameters that directly impacts the convergence of CDA. It is shown in \cite{nesterov2009primal} that an $O(\frac{1}{\sqrt{t}})$ convergence rate is ensured when $a_t$ decreases at $O(\frac{1}{\sqrt{t}})$ for nonsmooth objective functions. When the objective is smooth, an appropriate constant $a_t=a$ can be used to achieve an $O(\frac{1}{t})$ rate \cite{lu2018relatively}.  

For the projection operator in \eqref{cda}, a standard result in convex analysis {(Lemma 1 in \cite{nesterov2009primal})} is recalled in the following lemma.

\begin{lemma} \label{gamma_continuity}
	For any $u,v\in\mathbb{R}^m$, we have
	\begin{equation*}
	\begin{split}
	&\big\lVert \arg \min_{x\in \mathcal{X}}\big\{ a_t \langle u,x \rangle+d(x) \big\} -\arg \min_{x\in \mathcal{X}}\big\{ a_t \langle  v,x \rangle+ d(x) \big\} \big\rVert \\
	& \leq a_t \lVert u-v \rVert.
	\end{split}
	\end{equation*}
\end{lemma}


\subsection{Design of A New Distributed Dual Averaging Scheme}
In the literature, several distributed dual averaging algorithms have been developed accounting for different communication patterns among agents. Generally speaking, they both involve iteratively estimating the global dual variable $\sum_{k=0}^{t} \triangledown f(x_{k})$ in \eqref{cda} in the following way:
\begin{equation*}
q_{i,t+1}=\sum_{j=1}^{n}p_{ij}q_{j,t}+\triangledown f_i(x_{i,t+1})
\end{equation*}
where $q_{i,t}$ is an estimate of $\sum_{k=0}^{t} \triangledown f(x_{k})$ locally maintained by agent $i$ at time $t$, and $x_{i,t}$ is local estimate about the global minimizer. However, it is shown in \cite{liu2019distributed} that $q_{i,t}$ does not necessarily converge to the dual variable. Therefore, the control sequence $\{a_t\}_{t\geq 0}$ has to be decreasing for a slow convergence rate, i.e., $O(\frac{1}{\sqrt{t}})$.

To possibly accelerate convergence using a constant control sequence, the global dual variable must be more accurately estimated. Motivated by this, we propose to track the global dual variable according to the following rule: 

\begin{subequations}\label{2nd_order_consensus}
	\begin{align}
	s_{i,t+1}&=\sum_{j=1}^{n}p_{ij}s_{j,t}+\triangledown f_i(x_{i,t+1})-\triangledown f_i(x_{i,t})\\
			h_{i,t+1}&=\sum_{j=1}^{n}p_{ij}h_{j,t}+ s_{i,t+1}-s_{i,t}.
	\end{align}
\end{subequations}
Note that when $h_{i,0}=s_{i,0}$ one gets the new dual estimate
\begin{equation}\label{variant_consensus}
\sum_{k=0}^{t+1}h_{i,k}=\sum_{j=1}^{n}p_{ij}\sum_{k=0}^{t}h_{j,k}+ s_{i,t+1}.
\end{equation}
Thanks to it, the estimate about the global minimizer can be generated as follows. 
\begin{equation}\label{primal}
	{x}_{i,t+1}= \arg \min_{x\in\mathcal{X}}\Big\{ \sum_{k=0}^{t}a\langle h_{i,k},x \rangle+d(x)\Big\}.
\end{equation}


Denote $
\overline{h}_t=\frac{1}{n}\sum_{i=1}^{n}{h}_{i,t}$, $
\overline{s}_t=\frac{1}{n}\sum_{i=1}^{n}{s}_{i,t}$, and $g(t)=\frac{1}{n}\sum_{i=1}^{n}\triangledown f_i(x_{i,t})$. The following conservation property holds true.
\begin{lemma} \label{conservation_property}
	If $h_{i,0}=s_{i,0}=\triangledown f_i(x_{i,0}),i\in\mathcal{V}$, then
	\begin{equation*}
		\overline{h}_{t+1}=\overline{s}_{t+1}=g_{t+1}.
	\end{equation*}
\end{lemma} 
\begin{pf}
	The proof follows from projecting \eqref{2nd_order_consensus} into the average space.
	\end{pf}

The proposed algorithm is summarized in the following.
%

Initialization: Set $t=0$, $x_{i,0}=\arg\min_{x\in\mathcal{X}}{d(x)}=0_m$, $h_{i,0}=s_{i,0}=\triangledown f_i(x_{i,0}),\forall i\in\mathcal{V}$. Choose a constant control sequence $a_t=a$.

Each agent $i\in\mathcal{V}$ (in parallel)

1)  Receives $s_{j,t}, h_{j,t},\forall j \in\mathcal{N}_i$;

2) Performs local computation in \eqref{2nd_order_consensus} and \eqref{primal};

3)  Broadcasts $s_{i,t+1}, h_{i,t+1}$ to $ j \in\mathcal{N}_i$;

4)  Sets $t=t+1$.

%
%
%
%
%
%
%

\section{Main Result}

First, we set up an auxiliary sequence $\{y_t\}_{t\geq 0}$ that evolves according to the following rule
\begin{equation}\label{auxiliary_sequence}
	y_{t+1}= \arg \min_{x\in\mathcal{X}}\Big\{ \sum_{k=0}^{t}a\langle {g}_{k},x \rangle+d(x)\Big\},
\end{equation}
where the initial vector $y_0=\arg\min_{x\in\mathcal{X}}{d(x)}=0_m$.
Then, the deviation between $\{x_{i,t}\}_{t\geq 0}$ and $\{y_t\}_{t\geq 0}$ is analyzed. Finally, the convergence of $\{y_t\}_{t\geq 0}$ to the global minimizer is shown.

Define 
\begin{equation*}
{\bf x}_t=\begin{bmatrix}
x_{1,t}\\ x_{2,t} \\ \vdots \\ x_{n,t}
\end{bmatrix}, {\bf h}_t=\begin{bmatrix}
h_{1,t}\\ h_{2,t} \\ \vdots \\ h_{n,t}
\end{bmatrix},{\bf s}_t=\begin{bmatrix}
s_{1,t}\\ s_{2,t} \\ \vdots \\ s_{n,t}
\end{bmatrix},{\bf \triangledown}_t=\begin{bmatrix}
\triangledown f_1(x_{1,t})\\ \triangledown f_2(x_{2,t}) \\ \vdots \\ \triangledown f_n(x_{n,t})
\end{bmatrix}.
\end{equation*}
and ${\bf z}_{k+1}=\sum_{l=0}^{k}a{\bf h}_l$. 

The following lemma establishes the relation between sequences $\{x_{i,t}\}_{t\geq 0}$ and $\{y_t\}_{t\geq 0}$; the deviation between them represents the consensus error to be  compensated in convergence rate analysis.

\begin{lem}\label{consensus_error_lemma}
	For 	
	\begin{equation*}
	E(a)=\begin{bmatrix}
	\beta & a \\
	L(\beta+1) & \beta+La
	\end{bmatrix},
		\end{equation*}
		where $\beta =\sigma_2(P)$,
	if $\rho(E(a))<1$, it holds that
	\begin{equation}\label{consensus_error}
	\begin{split}
	\sum_{k=0}^{t-1}\lVert{{\bf x}}_{k}-\mathbf{1}y_{k}\rVert^2\leq \frac{{n}}{\big(1-\rho(E(a))\big)^2}\sum_{j=0}^{t-1}\lVert  y_{j+1}-y_{j} \rVert^2.
	\end{split}
	\end{equation}
\end{lem}

\begin{pf}
Please refer to Appendix A.

\end{pf}

The following lemma plays a similar role with the well-known dual averaging inequality (Theorem 2 in \cite{nesterov2009primal}) for nonsmooth optimization in convergence analysis. However, it further makes use of the smoothness of the objective in order to provide a much tighter bound for a faster convergence rate.

\begin{lem}\label{inexact_dual_averaging_inter}
	For $\{y_t\}_{t\geq 0}$ generated by \eqref{auxiliary_sequence}, it holds
	\begin{equation}\label{inexact}
	\begin{split}
	\sum_{k=0}^{t-1}\langle a g_k,y_{k+1}-x^* \rangle\leq d(x^*)-\sum_{k=0}^{t-1}D_{d}(y_{k+1}-y_{k}).
	\end{split}
	\end{equation}
\end{lem}   
\begin{pf}
The proof is postponed to Appendix B.
\end{pf}

%

We are now in a position to present the main result.
\begin{thm}\label{inexact_thm}
	If $\rho(E(a))<1$ and 
	\begin{equation*}
	 {aL} + \frac{aL}{\big(1-\rho(E(a))\big)^2} \leq \frac{1}{2},
	\end{equation*} then 
\begin{equation}
\begin{split}
f(\tilde{y}_{t})  -f(x^*)\leq  \frac{nd(x^*)}{at},
\end{split}
\end{equation}
where $\tilde{y}_t=\frac{1}{t}\sum_{k=0}^{t-1}y_{k+1}$.
\end{thm}  
\begin{pf}    
Consider

\begin{equation}
\begin{split}
&\sum_{j=1}^{n}a\Big( f_j(y_{k+1}) -f_j(x^*)\Big)\\
\leq& \sum_{j=1}^{n} a \Big( \frac{L}{2}\lVert y_{k+1}-x_{j,k} \rVert^2\\
&+f_j(x_{j,k})+\langle \triangledown f_j(x_{j,k}),y_{k+1}-x_{j,k} \rangle- f_j(x^*) \Big)\\
\leq&  \sum_{j=1}^{n}  a\Big(\frac{L}{2}\lVert y_{k+1}-x_{j,k} \rVert^2 +\langle \triangledown f_j(x_{j,k}),y_{k+1}-x^* \rangle \Big)  \\
=& \sum_{j=1}^{n}a\Big(\frac{L}{2}\lVert y_{k+1}-x_{j,k} \rVert^2 \Big)+ n  \langle  ag_{k}, {y}_{k+1}-x^*\rangle \\
\leq &   \frac{a}{2}{L}\lVert \mathbf{1}y_{k+1}-{\bf x}_{k} \rVert^2+n  \langle  ag_{k}, {y}_{k+1}-x^*\rangle\\
\leq &   \frac{a}{2}{L}\lVert \mathbf{1}y_{k+1}-\mathbf{1}y_{k}+\mathbf{1}y_{k}-{\bf x}_{k} \rVert^2+n  \langle  ag_{k}, {y}_{k+1}-x^*\rangle \\
\leq &  {a}L\Big( n\lVert y_{k+1}-y_{k} \rVert^2+\lVert \mathbf{1}y_{k}-{\bf x}_{k} \rVert^2 \Big)+n  \langle  ag_{k}, {y}_{k+1}-x^*\rangle,
\end{split}
\end{equation}
where the first inequality follows from the use of Lipschitz continuity of the gradient.

This together with convexity of $f_j$ allows us to further get
\begin{equation*}
\begin{split}
&at\Big(f(\tilde{y}_{t})  -f(x^*)\Big)
\leq at\Big(\sum_{j=1}^{n} f_j(\tilde{y}_{t})  -f(x^*)\Big)\\
\leq& \sum_{k=0}^{t-1}\sum_{j=1}^{n}a\Big( f_j(y_{k+1}) - f_j(x^*) \Big) \\
\leq &  {aLn}\sum_{k=0}^{t-1}\lVert y_{k+1}-y_{k} \rVert^2+n  \sum_{k=0}^{t-1}\langle  ag_{k}, {y}_{k+1}-x^*\rangle\\
&+{aL}\sum_{k=0}^{t-1}\lVert \mathbf{1}y_{k}-{\bf x}_{k} \rVert^2 \\
\leq & \Big(  {aLn} + \frac{aL{n}}{\big(1-\rho(E(a))\big)^2} \Big)\sum_{k=0}^{t-1}\lVert y_{k+1}-y_{k} \rVert^2\\
&+nd(x^*)-n\sum_{k=0}^{t-1}D_{d}(y_{k+1}-y_{k}).
\end{split}
\end{equation*}
Due to
$
D_d(y_{k+1},y_k)\geq \frac{1}{2} \lVert y_{k+1}-y_{k} \rVert^2,
$
we arrive at
\begin{equation}\label{rate_part1}
\begin{split}
&at\Big(f(\tilde{y}_{t})  -f(x^*)\Big)\\
\leq & {n}\Big(  {aL} + \frac{aL}{\big(1-\rho(E(a))\big)^2} -\frac{1}{2}\Big)\sum_{k=0}^{t-1}\lVert y_{k+1}-y_{k} \rVert^2+nd(x^*),
\end{split}
\end{equation}
thereby completing the proof.
\end{pf}
\begin{remark}
	Theorem \ref{inexact_thm} states that $\tilde{y}_t$ converges to the global minimizer at an $O(\frac{1}{t})$ rate. By \eqref{consensus_error} and convexity of the $2$-norm operator, one has
		\begin{equation}\label{consensus_result}
	\begin{split}
&	{t}\lVert{\tilde{\bf x}}_{t}-\mathbf{1}\tilde{y}_{t}\rVert^2\leq	\sum_{k=1}^{t}\lVert{{\bf x}}_{k}-\mathbf{1}y_{k}\rVert^2\\
	&\leq \frac{{n}}{\big(1-\rho(E(a))\big)^2}\sum_{j=0}^{t}\lVert  y_{j+1}-y_{j} \rVert^2,
	\end{split}
	\end{equation}
	where $\tilde{\bf x}_{t}=\frac{1}{t}\sum_{k=0}^{t-1}{\bf x}_{k+1}$. Moreover, from \eqref{rate_part1}, we know that the right-hand side of \eqref{consensus_result} remains finite as $t$ approaches infinity. Therefore, $\lVert{\tilde{ x}}_{i,t}-\tilde{y}_{t}\rVert^2$ converges at an $O(\frac{1}{t})$ rate, where $\tilde{ x}_{i,t}=\frac{1}{t}\sum_{k=0}^{t-1}{ x}_{i,k+1}$. This implies that $\tilde{ x}_{i,t}$ shares a similar convergence guarantee with $\tilde{ y}_{t}$.
	\end{remark}

\section{Simulation}

To verify the proposed method, we apply it to a large-scale LASSO problem. In this problem, the data tuple $(y_i,A_i)$ available at each agent $i\in\mathcal{V}$ satisfies the following equation:
 \begin{equation*}
 y_i=A_ix^*+b_i,
 \end{equation*} 
 where $A_i\in \mathbb{R}^{p_i\times m}$, $y_i\in\mathbb{R}^{p_i} $ , and $b_i\in\mathbb{R}^{m}$ is the additive Gaussian noise with zero mean and variance $\sigma^2$. 
 Usually, $p_i\ll m$ and $x^*$ is sparse. To recover $x^*$,
the following distributed optimization problem is considered:
\begin{equation*}
\min_{x\in\mathbb{R}^m}  \sum_{i=1}^{n}\frac{1}{2}\|y_i-A_ix \|^2, \quad \mathrm{s.t.} \,\, \lVert x \rVert_1 \leq R.
\end{equation*}

In the simulation, we set $n=50$, $m=10000$, $p_i=20, \forall i\in\mathcal{V}$. The matrix $A_i$ is randomly generated with $\mathcal{N}(0,1)$ elements. The minimizer $x^*$ is a sparse vector that only has $50$ non-zero $\mathcal{N}(0,1)$ entries. 
The variance for noise $b_i,\forall i\in\mathcal{V}$ is set as $\sigma^2=0.01$. Set $R= 1.1*\lVert x^* \rVert_1$. The communication network is characterized by an Erdos-Renyi graph with a $0.1$ connectivity ratio, and the doubly stochastic matrix $P$ associated with the graph is derived by following the Metropolis-Hastings rule.

For the purpose of comparison, the distributed projected gradient method (DPG) in \cite{nedic2010constrained}, and the distributed dual averaging (DDA) in \cite{duchi2011dual}
are simulated. To accommodate the theoretical results developed therein, the stepsize for DPG is chosen as $\frac{1}{\sqrt{t}}$; the control sequence in DDA is set as $a_t = \frac{1}{\sqrt{t}}$.
The control sequence for the proposed new DDA (N-DDA) is set as $a=\frac{1}{m}$. The initial primal variable for DPG 
is set as $x_{i,0}=0, \forall i \in\mathcal{V}$.

The simulation results are reported in the following. The performance is evaluated in terms of two criteria, that is, the primal variable residual of the first agent, i.e., $\frac{\lVert{x}_{1,t}-{x}^*\rVert^2}{\lVert{x}^*\rVert^2}$, and the objective value over the number of local iteration times. 
The results suggest that the proposed N-DDA enjoys a faster convergence rate. This is compatible with the theoretical results that PGA and DDA have a rate of $O(\frac{1}{\sqrt{t}})$ while N-DDA converges at an $O(\frac{1}{{t}})$ rate.

\begin{figure}
	\begin{center}
		\includegraphics[width=8.4cm]{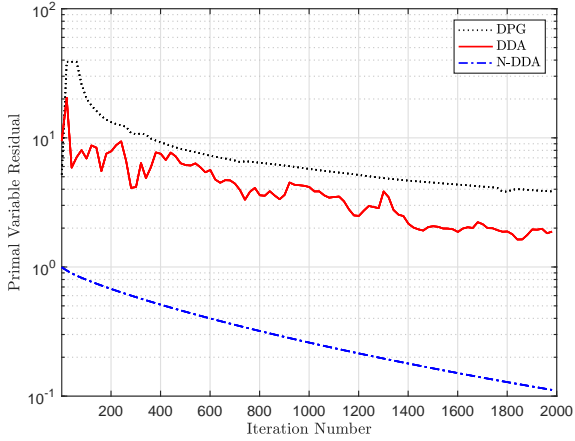}    
		\caption{Convergence of the primal variable residual.} 
		\label{fig:primal}
	\end{center}
\end{figure}

\begin{figure}
\begin{center}
\includegraphics[width=8.4cm]{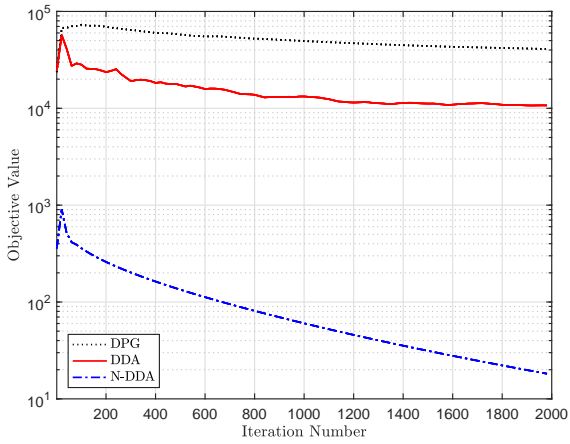}    
\caption{Convergence of the objective value.} 
\label{fig:objective}
\end{center}
\end{figure}

\section{Conclusion}
In this work, we proposed a new distributed dual averaging method tailored for smooth problems that has a convergence rate of $O(\frac{1}{t})$. This is made possible by a second-order consensus scheme that provides an accurate local estimate of the dual variable and a new analysis framework for dual averaging with inexact gradients. This work opens several avenues for future research, including the extension to smooth and strongly convex problems, and dynamic communication networks.


\bibliography{ifacconf}             

\appendix
\section{Proof of Lemma \ref{consensus_error_lemma}}    
	Since ${\bf s}_0={\bf h}_0=\triangledown_0$, from \eqref{2nd_order_consensus} we have
\begin{equation*}
\begin{split}
&{\bf z}_{k} =P{\bf z}_{k-1}+a{\bf s}_{k-1}.
\end{split}
\end{equation*}
By subtracting $\sum_{l=0}^{k-1}ag_l$ on both sides and the triangle inequality, we get
\begin{equation}\label{dual_error}
\begin{split}
&\lVert {\bf z}_k-\mathbf{1}\sum_{l=0}^{k-1}ag_l\rVert\\
\leq&  \lVert P{\bf z}_{k-1}-\mathbf{1}\sum_{l=0}^{k-2}ag_l\rVert +a\lVert {\bf s}_{k-1}-\mathbf{1}g_{k-1}\rVert \\
\leq& \beta 
\lVert {\bf z}_{k-1}-\mathbf{1}\sum_{l=0}^{k-2}ag_l\rVert +a\lVert {\bf s}_{k-1}-\mathbf{1}g_{k-1}\rVert.
\end{split}
\end{equation}
Similarly, it holds that
\begin{equation}
\begin{split}
&\lVert{{\bf s}}_k-\mathbf{1}g_k\rVert\\
=&\lVert P{{\bf s}}_{k-1}-\mathbf{1}g_{k-1}+\triangledown_{k}- \triangledown_{k-1}-\mathbf{1}g_k+\mathbf{1}g_{k-1}\rVert \\
\leq&  \beta \lVert{{\bf s}}_{k-1}-\mathbf{1}g_{k-1}\rVert + L\lVert {{\bf x}}_{k}-{\bf x}_{k-1}  \rVert,
\end{split}
\end{equation}
where the fact
\begin{equation*}
\lVert \triangledown_{k}- \triangledown_{k-1}-\mathbf{1}g_k+\mathbf{1}g_{k-1}\rVert \leq    \lVert \triangledown_{k}- \triangledown_{k-1}  \rVert
\end{equation*} 
and the Lipschitz continuity of the gradient are used to get the last inequality. Using Lemma \ref{gamma_continuity} over $\lVert {\bf x}_k-\mathbf{1}y_k\rVert$ and $\lVert {\bf x}_{k-1}-\mathbf{1}y_{k-1}\rVert$, and \eqref{dual_error} allows us to further get

\begin{equation}\label{gradient_error}
\begin{split}
&\lVert{{\bf s}}_k-\mathbf{1}g_k\rVert\\
\leq & \beta \lVert{{\bf s}}_{k-1}-\mathbf{1}g_{k-1}\rVert + L\lVert {\bf x}_k-\mathbf{1}y_k\rVert+ L\lVert {\bf x}_{k-1}-\mathbf{1}y_{k-1}\rVert \\
&+ L\lVert \mathbf{1}y_k-\mathbf{1}y_{k-1}\rVert\\
\leq &(\beta+La)\lVert{{\bf s}}_{k-1}-\mathbf{1}g_{k-1}\rVert + \sqrt{n}L\lVert y_k-y_{k-1}\rVert \\
&+(L+\beta L)\lVert {\bf z}_{k-1}-\mathbf{1}\sum_{l=0}^{k-2}ag_l\rVert.
\end{split}
\end{equation}

From \eqref{dual_error} and \eqref{gradient_error}, the following linear system inequality can be established:
\begin{equation}
\begin{split}
\begin{bmatrix}
\lVert {\bf z}_k-\mathbf{1}\sum_{l=0}^{k-1}ag_l\rVert \\
\lVert{{\bf s}}_k-\mathbf{1}g_k\rVert
\end{bmatrix} \leq& E(a) \begin{bmatrix}
\lVert {\bf z}_{k-1}-\mathbf{1}\sum_{l=0}^{k-2}ag_l\rVert \\
\lVert{{\bf s}}_{k-1}-\mathbf{1}g_{k-1}\rVert
\end{bmatrix} \\
& + \sqrt{n}L \begin{bmatrix}
0\\
\lVert y_k-y_{k-1}\rVert 
\end{bmatrix}.
\end{split}
\end{equation}

Since ${\bf h}_0={\bf s}_0=\triangledown_0=\mathbf{1}g_0$ by initialization, it holds that
\begin{equation}
\begin{split}
\begin{bmatrix}
\lVert {\bf z}_k-\mathbf{1}\sum_{l=0}^{k-1}ag_l\rVert \\
\lVert{{\bf s}}_k-\mathbf{1}g_k\rVert
\end{bmatrix} \leq& \sqrt{n}L \sum_{j=0}^{k-1}(E(a))^{k-j-1}\begin{bmatrix}
0\\
\lVert y_{j+1}-y_{j}\rVert 
\end{bmatrix} .
\end{split}
\end{equation}

It is easy to check that the eigenvalues of $E(a)$ are
\begin{equation*}
\frac{2\beta +aL \pm \sqrt{a^2L^2+4(\beta+1)aL}}{2}.
\end{equation*}
Since $\rho(E(a))<1$, one readily has $aL<\beta+1$.
Then, according to \cite{williams1992n},
\begin{equation}\label{consensus_error_mid}
\begin{split}
&\lVert {\bf z}_k-\mathbf{1}\sum_{l=0}^{k-1}ag_l\rVert\\
& \leq
\frac{\sqrt{n}aL}{\sqrt{a^2L^2+4(\beta+1)aL}}\sum_{j=0}^{k-1}(\lambda_1^{k-j-1}-\lambda_2^{k-j-1})\lVert y_{j+1}-y_{j}\rVert\\
& \leq  \sqrt{n}\sum_{j=0}^{k-1}\rho(E(a))^{k-j-1}\lVert y_{j+1}-y_{j}\rVert,
\end{split}
\end{equation}
where $\lambda_1>\lambda_2$ are eigenvalues of $E(a)$.
Therefore 
\begin{equation*}
\begin{split}
&\sum_{k=0}^{t-1}\lVert {\bf z}_k-\mathbf{1}\sum_{l=0}^{k-1}ag_l\rVert^2 \\
\leq&  {n}\sum_{k=1}^{t-1}\Big(\sum_{j=0}^{k-1}\rho(E(a))^{k-j-1}
\lVert y_{j+1}-y_{j}\rVert\Big)^2 \\
\leq & {n}\sum_{k=1}^{t-1}\\
&\Big(\sum_{j=0}^{k-1}\Big(\rho(E(a))^{\frac{k-j-1}{2}}\Big)^2\sum_{j=0}^{k-1}\Big(\rho(E(a))^{\frac{k-j-1}{2}}\lVert y_{j+1}-y_{j}\rVert\Big)^2 \\
\leq & {n}\sum_{k=1}^{t-1}\frac{1}{1-\rho(E(a))}\sum_{j=0}^{k-1}\rho(E(a))^{{k-j-1}}\lVert  y_{j+1}-y_{j}\rVert^2 \\
= &\frac{{n}}{1-\rho(E(a))}\sum_{k=1}^{t-1}\sum_{j=0}^{k-1}\rho(E(a))^{{k-j-1}}\lVert  y_{j+1}-y_{j} \rVert^2 \\
\leq& \frac{{n}}{\big(1-\rho(E(a))\big)^2}\sum_{j=0}^{t-1}\lVert  y_{j+1}-y_{j} \rVert^2,
\end{split}
\end{equation*}
which together with Lemma \ref{gamma_continuity} yields \eqref{consensus_error}. 
\section{Proof of Lemma \ref{inexact_dual_averaging_inter}}              
Define 
\begin{equation*}
m_{k}({x})= \langle \sum_{l=0}^{k}ag_l,x \rangle +d({x}) .
\end{equation*}
We then have
\begin{equation*}
m_k({x}) = m_{k-1}({x})+\langle ag_k,x \rangle
\end{equation*}
According to the definition of Bregman divergence, we have
\begin{equation*}
\begin{split}
&D_{m_{k-1}}(y_{k+1}-y_{k})\\
=&m_{k-1}(y_{k+1})-m_{k-1}(y_{k})-\langle \triangledown m_{k-1}(y_{k}), y_{k+1}-y_{k}\rangle
\end{split}
\end{equation*}
which is equivalent to
\begin{equation*}
\begin{split}
&D_{d}(y_{k+1}-y_{k})\\
=&m_{k-1}(y_{k+1})-m_{k-1}(y_{k})-\langle \triangledown m_{k-1}(y_{k}), y_{k+1}-y_{k}\rangle.
\end{split}
\end{equation*}
Since 
\begin{equation*}
y_{k}=\arg \min_{x\in\mathcal{X}} m_{k-1}(x),
\end{equation*}
by the optimality condition we have
\begin{equation*}
\langle \triangledown m_{k-1}(y_{k}),y_{k+1}-y_{k}\rangle \geq 0
\end{equation*}
and therefore
\begin{equation*}
\begin{split}
0&\leq m_{k-1}(y_{k+1})-m_{k-1}(y_{k})-D_{d}(y_{k+1}-y_{k}) \\
&= m_k(y_{k+1})-\langle ag_k,y_{k+1} \rangle -m_{k-1}(y_{k})-D_{d}(y_{k+1}-y_{k}) 
\end{split}
\end{equation*}
which is equivalent to
\begin{equation*}
\begin{split}
\langle a g_k,y_{k+1} \rangle\leq  m_k(y_{k+1}) -m_{k-1}(y_{k})-D_{d}(y_{k+1}-y_{k}) .
\end{split}
\end{equation*}
Summing the above equation over $k$ from $0$ to $t-1$ yields
\begin{equation}\label{inexact_gradient_result1}
\begin{split}
&\sum_{k=0}^{t-1}\langle a g_k,y_{k+1} \rangle\\
\leq&  m_{t-1}(y_{t}) -m_{0}(y_{1})+\langle a g_0,y_{1} \rangle-\sum_{k=0}^{t-1}D_{d}(y_{k+1}-y_{k}) \\
=& m_{t-1}(y_{t}) -\sum_{k=0}^{t-1}D_{d}(y_{k+1}-y_{k})
\end{split}
\end{equation}
We turn to consider
\begin{equation*}
\begin{split}
\sum_{k=0}^{t-1}\langle a g_k,-x^* \rangle&\leq  \max_{{x}\in\mathcal{X}} \Big\{ \langle \sum_{k=0}^{t-1} a g_k,-x \rangle-d(x)\Big\}+d(x^*) \\
& = -\min_{{x}\in\mathcal{X}} \Big\{ \langle \sum_{k=0}^{t-1} a g_k,x \rangle+d(x)\Big\}+d(x^*) \\
& = -m_{t-1}(y_{t})+d(x^*),
\end{split}
\end{equation*}
which in conjunction with \eqref{inexact_gradient_result1} gives rise to the inequality in \eqref{inexact}.
\end{document}